\documentclass{amsart}
\usepackage{amssymb}
\usepackage{amsmath}
\usepackage{a4}
\begin{document}
%\date{\version}
\newtheorem{theorem}{Theorem}[section]
\newtheorem{lemma}[theorem]{Lemma}
\newtheorem{remark}[theorem]{Remark}
\newtheorem{definition}[theorem]{Definition}
\newtheorem{corollary}[theorem]{Corollary}
\newtheorem{example}[theorem]{Example}
\def\qedbox{\hbox{$\rlap{$\sqcap$}\sqcup$}}
\makeatletter
  \renewcommand{\theequation}{%
   \thesection.\alph{equation}}
  \@addtoreset{equation}{section}
 \makeatother
%\date{\version}
%\def\BB{\mathcal{B}}
\def\MM{\mathfrak{M}}
\title[Manifolds with commuting Curvature Operators]
{The global geometry of Riemannian manifolds with commuting curvature operators}
\author[Brozos-V\'azquez]{M. Brozos-V\'azquez}
\address{Department of Geometry and Topology, Faculty of Mathematics, University of
Santiago de Compostela, 15782 Santiago de Compostela, Spain}
\email{mbrozos@usc.es}
\author[Gilkey]{P. Gilkey}
\address{Mathematics Department, University of Oregon, Eugene, OR 97403, USA}
\email{gilkey@.uoregon.edu}
\begin{abstract} We give manifolds whose Riemann curvature operators commute, i.e. which satisfy
$\mathcal{R}(x_1,x_2)\mathcal{R}(x_3,x_4)=\mathcal{R}(x_3,x_4)\mathcal{R}(x_1,x_2)$
for all tangent vectors $x_i$ in both the
Riemannian and the higher signature settings. These manifolds have global geometric phenomena which are quite different for
higher signature manifolds than they are for Riemannian manifolds.
Our focus is on global properties; questions of geodesic
completeness and the behaviour of the exponential map are
investigated.
\end{abstract}
\maketitle
\section{Introduction}

The study of curvature is central not only to Differential Geometry but to Global Analysis and Topology as well;
one must relate properties of the curvature tensor to the underlying geometry and topology of the manifold
under consideration -- the curvature tensor not only describes the geometry, but in a broad range of situations provides
useful information about the topology and analytic properties of the manifold, especially when considering the global aspects.
Much of this analysis involves examining a natural operator associated with the curvature; one studies commutativity or
spectral properties of this operator. There are many operators which can be considered; the skew-symmetric curvature operator
$\mathcal{R}(x,y)$ is perhaps the most natural one in this context. In this paper, we investigate a very natural geometric
question by studying the global properties of manifolds which satisfy the condition
$$\mathcal{R}(x_1,x_2)\mathcal{R}(x_3,x_4)=\mathcal{R}(x_3,x_4)\mathcal{R}(x_1,x_2)$$
for all tangent vectors $x_i$ -- i.e. manifolds where the curvature operator is totally commutative.

The subject can properly have been said to have started with
Osserman \cite{O90}; seminal papers by Bla\v zi\' c et. al.
\cite{BBR00}, by  Bonome et al. \cite{refBCG}, by Chi
\cite{Chin88}, by Ivanov and Petrova \cite{IP98}, by Ivanova and
Stanilov \cite{IS95}, by  Nikolayevsky \cite{Ni04c}, by Stanilov
\cite{S04}; by Stanilov and Videv \cite{SV92}, and by Tsankov
\cite{Y05} are only a few that could be mentioned; as the
literature is a vast one, we must content ourselves for referring
to the bibliographies in \cite{GF02,Gl02} for further citations --
the field is a vibrant one with a growing bibliography! A
variety of methods, ranging from algebraic topology to elliptic
operator theory and to classical invariance theory, have
been used to study the spectral geometry of the curvature tensor.

Although the focus of this paper is on questions in global geometry, it is often convenient to study geometrical problems by first working in
a purely algebraic context. One says that
$\mathfrak{M}:=(V,\langle\cdot,\cdot\rangle,A)$ is a
$0$-model if $V$ is a vector space of dimension $m$, if $\langle\cdot,\cdot\rangle$ is a non-degenerate symmetric bilinear form on
$V$ of signature $(p,q)$, and if
$A\in\otimes^4V^*$ is an algebraic curvature tensor on $V$; this means that $A$ satisfies the usual curvature symmetries:
\begin{eqnarray*}
&&A(x,y,z,w)=A(z,w,x,y)=-A(y,x,z,w),\\
&&0=A(x,y,z,w)+A(y,z,x,w)+A(z,x,y,w)\,.
\end{eqnarray*}
 Let $\mathcal{A}$ be the associated
skew-symmetric curvature operator:
$$\langle\mathcal{A}(x,y)z,w\rangle=A(x,y,z,w)\,.$$
One says that $A$ is {\it skew Tsankov} if
\begin{equation}\label{eqn-1.a}
\mathcal{A}(x,y)\mathcal{A}(u,v)=\mathcal{A}(u,v)\mathcal{A}(x,y)\quad
\forall\; x,y,u,v\,.
\end{equation}

If $\mathcal{M}:=(M,g)$ is a pseudo-Riemannian manifold, let $R$ be the curvature tensor of the Levi-Civita
connection and let the associated
$0$-model be given by:
$$\MM(\mathcal{M},P):=(T_PM,g_P,R_P)\,.$$
One says that
$\mathcal{M}$ is {\it skew Tsankov}  if $\MM(\mathcal{M},P)$ is skew Tsankov for every $P\in M$. The notation is motivated by the seminal
result of Tsankov \cite{Y05} who, following foundational suggestions of Stanilov, studied similar questions for hypersurfaces in
$\mathbb{R}^m$; Tsankov imposed an extra condition of orthogonality that we shall not impose here. Similar questions arise for the Jacobi
operator; see
\cite{BG05} for further details.

Here is a brief guide to this paper. In Section
\ref{sect-2}, we give a complete classification of Riemannian ($p=0$) skew Tsankov algebraic curvature tensors. In Section \ref{sect-3},
we present a family of irreducible $3$-dimensional Riemannian skew Tsankov manifolds. In Section \ref{sect-4}, we present a family of
irreducible $4$-dimensional Riemannian skew Tsankov manifolds. These examples indicate that despite the fact that the algebraic
classification is complete, the geometric classification promises to be more difficult; questions of global geometry turn out to
be deeper than the corresponding algebraic questions in this instance.

We then turn our attention to the pseudo-Riemannian setting. It
turns out that many manifolds which appeared in other settings are
also skew Tsankov. Dunn manifolds were introduced in \cite{DG05};
they are a family of neutral signature pseudo-Riemannian manifolds
which are Osserman, Ivanov--Petrova, and Szab\'o. In Section
\ref{sect-5}, we show these manifolds are skew Tsankov with
skew-symmetric curvature operators which are nilpotent of order
$2$. Certain Fiedler manifolds have been shown to be nilpotent
Osserman of arbitrarily high order \cite{FG03}. In Section
\ref{sect-6}, we show that Fiedler manifolds are also skew Tsankov
and that their skew-symmetric curvature operator is nilpotent of
order $3$. Certain  Nik\v cevi\'c manifolds are skew Tsankov as
well; the verification that these manifolds in \cite{GN04a} are
skew Tsankov in dimension $m=6$ and signature $(4,2)$ is
relatively straightforward and will be omitted in the
interests of brevity. Thus there are many examples of skew
Tsankov manifolds in the higher signature context; these examples
indicate that even in the algebraic setting, the classification is
likely to be far more complicated and this is a fruitful subject
for further inquiry.

If $V=V_1\oplus V_2$ is a non-trivial orthogonal direct sum
decomposition of $V$ which induces a decomposition $A=A_1\oplus
A_2$, then $\MM=(V,\langle\cdot,\cdot\rangle,A)$ is said to be
decomposable and we write $\MM=\MM_1\oplus\MM_2$ where
$\MM_i:=(V_i,\langle\cdot,\cdot\rangle|_{V_i},A_i)$; $\MM$
is said to be indecomposable otherwise. A pseudo-Riemannian
manifold $\mathcal{M}=(M,g)$ is said to be reducible at a point
$P\in M$ if there is a neighborhood $\mathcal{O}$ of $P$ in $M$
and a Cartesian product
$\mathcal{O}=\mathcal{O}_1\times\mathcal{O}_2$ which induces an
orthonormal decomposition $g_{\mathcal{O}}=g_{\mathcal{O}_1}\oplus
g_{\mathcal{O}_2}$; $\mathcal{M}$ is locally
irreducible at $P$ if $\mathcal{M}$ is not reducible at $P$.

The manifolds of Section \ref{sect-5} and \ref{sect-6} are all geodesically complete; by contrast, the scalar curvature of the manifolds
described in Sections \ref{sect-3} and \ref{sect-4} blows up in finite time along certain geodesics and thus these manifolds are necessarily
geodesically incomplete and can not be embedded isometrically in a geodesically complete manifold. It is not known if there are any
irreducible skew Tsankov Riemannian manifolds of dimension $m\ge3$ which are geodesically complete.

\section{The classification of Riemannian skew Tsankov algebraic curvature tensors}\label{sect-2}
\begin{theorem}\label{thm-2.1}
Let $\MM:=(V,\langle\cdot,\cdot\rangle,A)$ be a Riemannian $0$-model.
\begin{enumerate}\item $\MM$ is
skew Tsankov if and only if there exists a orthogonal
direct sum decomposition $V=V_1\oplus...\oplus V_k\oplus W$
decomposing $A=A_1\oplus...\oplus A_k\oplus 0$ where $\dim(V_i)=2$
for all i. \item  $\MM$ is skew Tsankov and indecomposable if and
only if $\dim(V)=2$ and $A\ne0$.
\end{enumerate}\end{theorem}

\begin{proof}Suppose given an orthogonal direct sum decomposition $V=V_1\oplus...\oplus V_k\oplus W$ so $A=A_1\oplus...\oplus
A_k\oplus 0$ where $\dim(V_i)=2$ for $1\le i\le k$. Let $\{e_i^1,e_i^2\}$ be an orthonormal basis for $V_i$. Given $x,y\in V$, there
exist coefficients
$\varepsilon_i(x,y)\in\mathbb{R}$ with
\begin{equation}\label{eqn-2.a}
\mathcal{A}(x,y)\xi=\left\{\begin{array}{rll}
-\varepsilon_i(x,y)e_i^2&\text{ if }&\xi=e_i^1,\\
 \varepsilon_i(x,y)e_i^1&\text{ if }&\xi=e_i^2,\\
0&\text{ if }&\xi\perp\operatorname{Span}\{e_i^1,e_i^2\}\,.
\end{array}\right.\end{equation}
We may then show $\MM$ is skew Tsankov by computing
$$
\mathcal{A}(x,y)\mathcal{A}(\bar x,\bar y)\xi=\left\{\begin{array}{rll}
-\varepsilon_i(x,y)\varepsilon_i(\bar x,\bar y)e_i^1&\text{ if }&\xi=e_i^1,\\
-\varepsilon_i(x,y)\varepsilon_i(\bar x,\bar y)e_i^2&\text{ if }&\xi=e_i^2,\\
0&\text{ if }&\xi\perp\operatorname{Span}_{1\le i\le k}\{e_i^1,e_i^2\}\,.
\end{array}\right.$$

Conversely, suppose that $\MM$ is skew Tsankov. One may simultaneously skew-diagonalize the collection
$\{A(x,y)\}_{x,y\in V}$ of commuting skew-adjoint linear operators to find an orthonormal set
$\{e_i^1,e_i^2\}$ and to find functions $\varepsilon_i(x,y)$ where $1\le i\le k$ so that Equation (\ref{eqn-2.a})
holds.
Extend this to a full orthonormal basis
$\mathcal{B}:=\{e_1^1,e_1^2,...,e_k^1,e_k^2,f_1,...,f_l\}$ for $V$. Then the only non-zero entries
in the curvature tensor relative to this base are $A(\cdot,\cdot,e_i^1,e_i^2)$ modulo the usual
$\mathbb{Z}_2$ symmetry. Interchanging the first 2 entries with the last 2 entries shows the
only non-zero curvatures are
$$A(e_i^1,e_i^2,e_j^1,e_j^2)\,.$$
On the other hand, if $i\ne j$, we can use the Bianchi identity to express
$$
A(e_i^1,e_i^2,e_j^1,e_j^2)=
-A(e_i^1,e_j^1,e_j^2,e_i^2)
-A(e_i^1,e_j^2,e_i^2,e_j^1)=0
$$
and thus the only non-zero curvatures are $a_i:=A(e_i^1,e_i^2,e_i^2,e_i^1)$. Thus setting $V_i:=\operatorname{Span}\{e_i^1,e_i^2\}$ yields
the desired decomposition $A=A_1\oplus...\oplus A_k\oplus 0$. Assertions (1) and (2) now follow.
\end{proof}

\section{$3$-dimensional irreducible skew Tsankov manifolds}\label{sect-3}

We construct irreducible $3$-dimensional examples by taking a product of the interval $(0,\infty)$ with a Riemann surface:

\begin{theorem}\label{thm-3.1} Let $\mathcal{N}:=(N,g_N)$ be a Riemann surface which does not have constant sectional curvature $+1$.
Give $\mathcal{M}:=((0,\infty)\times N,g_M)$ the warped product metric $g_M:=dt^2+t^2g_N$ for $t\in(0,\infty)$. Then $\mathcal{M}$ is an
irreducible skew Tsankov manifold with scalar curvature $\tau_{\mathcal{M}}=t^{-2}\{\tau_{\mathcal{N}}-2\}$.
\end{theorem}
\begin{proof} Choose isothermal coordinates to express
$ds^2_N=e^{2\alpha}(dx_1^2+dx_2^2)$, at least locally. Let
$\partial_1:=\partial_{x_1}$, let $\partial_2:=\partial_{x_2}$,
and let $\partial_3:=\partial_t$. Let
$\alpha_i:=\partial_i(\alpha)$ and
$\alpha_{ij}:=\partial_i\partial_j(\alpha)$. We have
$$g_M(\partial_1,\partial_1)=g_M(\partial_2,\partial_2)=t^2e^{2\alpha}\quad\text{and}\quad
g_M(\partial_3,\partial_3)=1\,.$$ The non-zero Christoffel
symbols of the first kind must have at least one repeated index
different from $3$:
$$\begin{array}{lll}
\Gamma_{111}=\alpha_1t^2e^{2\alpha},&
\Gamma_{112}=-\alpha_2t^2e^{2\alpha},&
\Gamma_{113}=-te^{2\alpha},\\
\Gamma_{121}=\Gamma_{211}=\alpha_2t^2e^{2\alpha},\quad&
\Gamma_{131}=\Gamma_{311}=te^{2\alpha},\quad\\
\Gamma_{221}=-\alpha_1t^2e^{2\alpha},&
\Gamma_{222}=\alpha_2t^2e^{2\alpha},&
\Gamma_{223}=-te^{2\alpha},\\
\Gamma_{122}=\Gamma_{212}=\alpha_1t^2e^{2\alpha},&
\Gamma_{322}=\Gamma_{232}=te^{2\alpha}\,.
\end{array}$$Since the metric is diagonal, we can raise indices to see:
\begin{eqnarray*}
&&\nabla_{\partial_1}\partial_1=\alpha_1\partial_1-\alpha_2\partial_2-te^{2\alpha}\partial_3,\\
&&\nabla_{\partial_1}\partial_2=\nabla_{\partial_2}\partial_1=\alpha_2\partial_1+\alpha_1\partial_2,\\
&&\nabla_{\partial_1}\partial_3=\nabla_{\partial_3}\partial_1=t^{-1}\partial_1,\\
&&\nabla_{\partial_2}\partial_2=-\alpha_1\partial_1+\alpha_2\partial_2-te^{2\alpha}\partial_3,\\
&&\nabla_{\partial_2}\partial_3=\nabla_{\partial_3}\partial_2=t^{-1}\partial_2\,.
\end{eqnarray*}
It is now an easy exercise to determine the curvature operator; we shall omit details in the interests of brevity. One has:
\begin{eqnarray*}
&&\mathcal{R}_{\mathcal{M}}(\partial_1,\partial_2)\partial_1=(\alpha_{11}+\alpha_{22}+e^{2\alpha})\partial_2,\\
&&\mathcal{R}_{\mathcal{M}}(\partial_1,\partial_2)\partial_2=-(\alpha_{11}+\alpha_{22}+e^{2\alpha})\partial_1,\\
&&\mathcal{R}_{\mathcal{M}}(\partial_1,\partial_3)\partial_3=\mathcal{R}_{\mathcal{M}}(\partial_2,\partial_3)\partial_3=0\,.
\end{eqnarray*}
As the only non-zero curvature is $R_{\mathcal{M}}(\partial_1,\partial_2,\partial_2,\partial_1)
=-t^2e^{2\alpha}(\alpha_{11}+\alpha_{22}+e^{2\alpha})$,
Theorem \ref{thm-2.1} implies that $\mathcal{M}$ is skew Tsankov. This calculation also yields
\begin{eqnarray*}
\tau_{\mathcal{M}}=t^{-2}\{-2e^{-2\alpha}(\alpha_{11}+\alpha_{22})-2\}\,.
\end{eqnarray*}

An analogous computation on $\mathcal{N}$ yields:
$$\begin{array}{lll}
\Gamma_{111}=\alpha_1t^2e^{2\alpha},&
\Gamma_{112}=-\alpha_2t^2e^{2\alpha},\quad&
\Gamma_{121}=\Gamma_{211}=\alpha_2t^2e^{2\alpha},\\
\Gamma_{221}=-\alpha_1t^2e^{2\alpha},\quad&
\Gamma_{222}=\alpha_2t^2e^{2\alpha},&\Gamma_{122}=\Gamma_{212}=\alpha_1t^2e^{2\alpha},
\end{array}$$
so the Christoffel symbols of the second kind and the curvature are given by:
$$\begin{array}{ll}
\nabla_{\partial_1}\partial_1=\alpha_1\partial_1-\alpha_2\partial_2,&
  \nabla_{\partial_2}\partial_2=-\alpha_1\partial_1+\alpha_2\partial_2,\\
\nabla_{\partial_1}\partial_2=\nabla_{\partial_2}\partial_1=\alpha_2\partial_1+\alpha_1\partial_2,&
  \mathcal{R}_{\mathcal{N}}(\partial_1,\partial_2)\partial_1=(\alpha_{11}+\alpha_{22})\partial_2\,.
\vphantom{\vrule height 11pt}\end{array}$$
Theorem \ref{thm-3.1} now follows; $\mathcal{M}$ is indecomposable because
$\operatorname{Range}\{\mathcal{R}\}=\operatorname{Span}\{\partial_1,\partial_2\}$ and because
$\tau_{\mathcal{M}}=t^{-2}(\tau_{\mathcal{N}}-2)$ exhibits non-trivial dependence on $t$.\end{proof}

\begin{remark}\label{rmk-3.2}
\ \begin{enumerate}\item \rm Let $f(x_1,x_2)$ be an isometric
embedding of a Riemann surface $N$ in $S^3\subset\mathbb{R}^4$.
Define an embedding of $(0,\infty)\times N$ in $\mathbb{R}^4$ by
setting $F(t,x):=tf(x)$. Theorem \ref{thm-3.1} may then be used to
see the resulting hypersurface in $\mathbb{R}^4$ is skew Tsankov;
such hypersurfaces appear in Tsankov \cite{Y05}. \item Choose a
point $x\in N$ where $\tau_{\mathcal{N}}(x)\ne 2$ and let
$\gamma_x(t):=t\times x$. Then $\gamma_x$ is a unit speed geodesic
and
$\lim_{t\rightarrow0}|\tau_{\mathcal{M}}(\gamma_x(t))|=\infty$.
Thus $\mathcal{M}$ exhibits scalar curvature blowup at finite
time. This shows $\mathcal{M}$ is geodesically incomplete and can
not be embedded isometrically in a geodesically complete manifold.
It is not known whether or not there exist irreducible complete
skew Tsankov $3$-dimensional manifolds.
\end{enumerate}\end{remark}

\section{$4$-dimensional irreducible Skew-Tsankov Manifolds}\label{sect-4}
We take a warped product metric with a flat base and a flat fiber. Denote
the usual coordinates on $\mathbb{R}^4$ by
$(x_1,x_2,x_3,x_4)$. Let $\partial_i:=\partial_{x_i}$ and let
$$\mathcal{O}:=\{(x_1,x_2,x_3,x_4)\in\mathbb{R}^4:x_3>0,x_4>0\}\,.$$

\begin{theorem}\label{thm-4.1}
For $\beta>0$,
let $\mathcal{M}_\beta:=(\mathcal{O},g_\beta)$ where
$$\begin{array}{ll}
g_\beta(\partial_1,\partial_1)=x_3^{2},\quad&
g_\beta(\partial_2,\partial_2)=(x_3+\beta x_4)^{2},\\
g_\beta(\partial_3,\partial_3)=1,&
g_\beta(\partial_4,\partial_4)=1\,.\vphantom{\vrule height 11pt}
\end{array}$$
\begin{enumerate}\item $\mathcal{M}_\beta$ is an indecomposable skew Tsankov manifold.
\item The scalar curvature $\tau_{\mathcal{M}_\beta}=-2x_3^{-1}(x_3+\beta x_4)^{-1}$.
\item $\mathcal{M}_{\beta_1}$ is not isometric to $\mathcal{M}_{\beta_2}$ for $\beta_1\ne\beta_2$.
\end{enumerate}\end{theorem}

\begin{proof} The non-zero Christoffel symbols are given by:
$$\begin{array}{ll}
\Gamma_{113}=-x_3,&\Gamma_{131}=\Gamma_{311}=x_3,\\
\Gamma_{223}=-(x_3+\beta x_4),&\Gamma_{232}=\Gamma_{322}=x_3+\beta x_4,\\
\Gamma_{224}=-\beta(x_3+\beta x_4),\quad&\Gamma_{242}=\Gamma_{422}=\beta(x_3+\beta x_4)\,.
\end{array}$$
Since the metric is diagonal, we may raise indices to compute:
$$\begin{array}{l}
\nabla_{\partial_1}\partial_1=-x_3\partial_3,\\
\nabla_{\partial_1}\partial_3=\nabla_{\partial_3}x_1=x_3^{-1}\partial_1,\\
\nabla_{\partial_2}\partial_2=-(x_3+\beta x_4)\partial_3-\beta(x_3+\beta x_4)\partial_4,\\
\nabla_{\partial_2}\partial_3=\nabla_{\partial_3}\partial_2=(x_3+\beta x_4)^{-1}\partial_2,\\
\nabla_{\partial_2}\partial_4=\nabla_{\partial_4}\partial_2=\beta(x_3+\beta x_4)^{-1}\partial_2\,.
\end{array}$$
The curvature operator can now be determined; as before, we shall omit the detailed computations in the interests of brevity:
\begin{eqnarray*}
&&\mathcal{R}(\partial_1,\partial_2)\partial_1=x_3(x_3+\beta x_4)^{-1}\partial_2,\\
&&\mathcal{R}(\partial_1,\partial_2)\partial_2=-x_3^{-1}(x_3+\beta x_4)\partial_1\,.
\end{eqnarray*}
The remaining curvatures vanish so the only non-zero curvature is
$$\mathcal{R}(\partial_1,\partial_2,\partial_2,\partial_1)=-x_3(x_3+\beta x_4)$$
and hence $\mathcal{M}$ is skew Tsankov by Theorem \ref{thm-2.1}. This establishes Assertion (1); Assertion (2) follows from the computations
performed above.

Let
$\mathcal{E}:=\operatorname{Range}\{\mathcal{R}\}=\operatorname{Span}\{\partial_1,\partial_2\}$
and let
$\mathcal{F}:=\mathcal{E}^\perp=\operatorname{Span}\{\partial_3,\partial_4\}$.
These spaces are invariantly defined. We have
\begin{eqnarray*}
&&\ln|\tau|=\ln(2)-\ln(x_3)-\ln(x_3+\beta x_4),\\
&&\nabla^2\{\ln|\tau|\}|_{\mathcal{F}}=\left(\begin{array}{ll}x_3^{-2}+(x_3+\beta x_4)^{-2}&\beta(x_3+\beta x_4)^{-2},\\
     \beta(x_3+\beta x_4)^{-2}&\beta^2(x_3+\beta x_4)^{-2}\end{array}\right),\\
&&\det(\nabla^2\{\ln|\tau|\}|_{\mathcal{F}})=\beta^2x_3^{-2}(x_3+\beta
x_4)^{-2}=\textstyle\frac14\beta\tau_{\mathcal{M}_\beta}^2\,.
\end{eqnarray*}
This shows that $\beta$ is an isometry invariant of $\mathcal{M}_\beta$. Furthermore since
$H|_\mathcal{F}$ has rank $2$, $\mathcal{M}$ is irreducible.
\end{proof}

\begin{remark}\label{rmk-4.2}
\rm As in the example described in Section \ref{sect-3}, the
scalar curvature blows up at finite time along the geodesic
$\gamma(t)=(1,1,t,1)$; thus $\mathcal{M}_\beta$ can not be
isometrically embedded as an open subset of a complete manifold;
it is not known whether or not every irreducible complete skew
Tsankov Riemannian manifold is necessarily
$2$-dimensional.\end{remark}

\section{Dunn manifolds}\label{sect-5}
We study the following family of {\it Dunn manifolds} which was first introduced in \cite{DG05} in a different context.
\begin{theorem}\label{thm-5.1} Let $(x_1,...,x_p,y_1,...,y_p)$ be
coordinates on $\mathbb{R}^{2p}$. Let
$\psi_{ij}(x)=\psi_{ji}(x)$ be given. Let
$\mathcal{M}:=(\mathbb{R}^{2p},g)$ be the manifold of neutral signature $(p,p)$ where
$$g(\partial_{x_i},\partial_{x_j})=\psi_{ij}(x)\quad\text{and}\quad
  g(\partial_{x_i},\partial_{y_i})=1\,.$$
Then $\mathcal{M}$ is
skew Tsankov and $\mathcal{R}$ is nilpotent of order $2$.
\end{theorem}

\begin{proof} The non-zero Christoffel symbols are:
\begin{eqnarray*}
&&g(\nabla_{\partial_{x_i}}\partial_{x_j},\partial_{x_k})=\textstyle\frac12
(\psi_{ik/j}+\psi_{jk/i}-\psi_{ij/k}),\\
&&\nabla_{\partial_{x_i}}\partial_{x_j}
   =\textstyle\frac12\sum_k(\psi_{ik/j}+\psi_{jk/i}-\psi_{ij/k})\partial_{y_k}\,.
\end{eqnarray*}
From this it follows that the possibly non-zero entries in
curvature tensor $R$ are:
$$R_{ijkl}
   =\textstyle-\frac12\sum_l(\psi_{il/jk}+\psi_{jk/il}-\psi_{ik/jl}-\psi_{jl/ik})\,.$$
Consequently
$\mathcal{R}(\partial_{x_i},\partial_{x_j})\partial_{x_k}=\textstyle\sum_\ell R_{ijkl}\partial_{y_\ell}$.
This shows that
\begin{eqnarray*}
&&\operatorname{Range}(\mathcal{R})\subset\operatorname{Span}\{\partial_{y_i}\}\quad\text{and}\quad
\operatorname{Span}\{\partial_{y_i}\}\subset\operatorname{Ker}(\mathcal{R})\,.
\end{eqnarray*}
Thus $\mathcal{R}(\xi_1,\xi_2)\mathcal{R}(\xi_3,\xi_4)=0$ for all $\xi_1,\xi_2,\xi_3,\xi_4$ so
$\mathcal{M}$ is skew Tsankov.\end{proof}

\begin{remark}\label{rmk-5.2}\rm
These manifolds have been studied extensively. These manifolds are all geodesically complete and the
exponential map is a global diffeomorphism. If $\psi_{ij}=\partial_{x_i}f\partial_{x_j}f$ for some function $f$, then $\mathcal{M}$ is
realizable as a hypersurface in $\mathbb{R}^{(p,p+1)}$. Certain manifolds in this family are curvature homogeneous but not homogeneous. We
refer to
\cite{DG05} for further details. Thus, in contrast to the Riemannian setting, there are examples which are global in the sense that
they are geodesically complete.\end{remark}

\section{Fiedler Manifolds}\label{sect-6}
The following family of examples was first introduced in
\cite{FG03}.

\begin{theorem}\label{thm-6.1}
Let $(x,u_1,...,u_\nu,y)$ be coordinates on $\mathbb{R}^{\nu+2}$.
Let $f\in C^\infty(\mathbb{R}^\nu)$ and let $\Xi=\Xi_{ab}$ be an invertible
symmetric
$\nu\times\nu$ matrix of signature $(r,s)$. Define a metric $g$ of signature $(r+1,s+1)$ on
$\mathbb{R}^{\nu+2}$ by setting:
$$\begin{array}{llll}
g(\partial_x,\partial_x)=-2f(\vec u),\quad&
g(\partial_x,\partial_y)=1,\quad&\text{and}&
g(\partial_{u_a},\partial_{u_b})=\Xi_{ab}.\vphantom{\vrule height 12pt}
\end{array}$$
Then $\mathcal{M}$ is skew Tsankov and $\mathcal{R}$ is nilpotent of order $3$.
\end{theorem}

\begin{proof}Since $d\Xi=0$, the potentially non-zero Christoffel symbols are:
\begin{eqnarray*}
&&g(\nabla_{\partial_x}\partial_x,\partial_{u_a})=\partial_{u_a}(f),\\
&&g(\nabla_{\partial_{u_a}}\partial_x,\partial_x)=
g(\nabla_{\partial_x}\partial_{u_a},\partial_x)=-\partial_{u_a}(f)\,.
\end{eqnarray*}
Let $\Xi^{ab}$ be the inverse matrix. Then
\begin{eqnarray*}
&&\nabla_{\partial_x}\partial_x=\textstyle\sum_{ab}\Xi^{ab}\partial_{u_a}(f)\partial_{u_b},
\\
&&\nabla_{\partial_x}\partial_{u_a}=\nabla_{\partial_{u_a}}\partial_x=-\partial_{u_a}(f)\partial_y\,.
\end{eqnarray*}
The quadratic terms
in the Christoffel symbols play no role in the calculation of $R$. Let $f_{ab}:=\partial_{u_a}\partial_{u_b}f$. The possibly
non-zero components of
$R$ and of $\mathcal{R}$ are
\begin{eqnarray*}
&&R(\partial_x,\partial_{u_a},\partial_{u_b},\partial_x)=f_{ab},\\
&&\mathcal{R}(\partial_x,\partial_{u_a})\partial_{u_b}=f_{ab}
\partial_y\quad\text{and}\quad
\mathcal{R}(\partial_x,\partial_{u_a})\partial_x=-\Xi^{bc}f_{ac}\partial_{u_b}\,.
\end{eqnarray*}
Thus the only potentially non-zero quadratic terms in the curvature are
$$\mathcal{R}(\partial_x,\partial_{u_d})\mathcal{R}(\partial_x,\partial_{u_a})\partial_x
=-\Xi^{bc}f_{ac}f_{db}\partial_y\,.$$
It now follows that
$\mathcal{R}(\partial_x,\partial_{u_d})\mathcal{R}(\partial_x,\partial_{u_a})\partial_x
  =\mathcal{R}(\partial_x,\partial_{u_a})\mathcal{R}(\partial_x,\partial_{u_d})\partial_x$.
This shows that $\mathcal{M}$ is skew-Tsankov and that $\mathcal{R}$ is nilpotent of order 3.
\end{proof}

These manifolds are irreducible
for generic $f$. They are complete for certain choices of the
warping function but are not in general geodesically complete. Of
particular interest is the Lorentzian case. We refer to
\cite{GN05} for the following results that again relate to the
global geometry of these examples. Let
$\mathcal{M}_{f}:=(\mathbb{R}^3,g_{f})$ where $g_{f}$ is the
Lorentz metric on $\mathbb{R}^3$ given by:
$$
g_f(\partial_x,\partial_x)=-2f(y)\quad\text{and}\quad g_f(\partial_x,\partial_{\tilde x})=g_f(\partial_y,\partial_y)=1\,.
$$

\begin{example}\label{exm-x}\rm
Let $\mathcal{S}_\varepsilon$ be defined by $f_\varepsilon(y):=\frac12\varepsilon y^2$ for $\varepsilon=\pm1$.
\begin{enumerate}
\item The manifolds $\mathcal{S}_\pm$ are geodesically complete.
\item The map $\exp_P$ for $\mathcal{S}_+$ is not surjective
$\forall P\in\mathbb{R}^3$.
\item The map $\exp_P$ for
$\mathcal{S}_-$ is a global diffeomorphism from $T_P\mathbb{R}^3$
to $\mathbb{R}^3$ $\forall P\in\mathbb{R}^3$.
\end{enumerate}\end{example}

We say that a pseudo-Riemannian manifold $\mathcal{M}$ is $k$-curvature
homogeneous if given any two points $P$ and $Q$ of $M$, there is
an isometry $\phi:T_PM\rightarrow T_QM$ so that
$\phi^*\nabla^iR_Q=\nabla^iR_P$ for $i\le k$. We say $\mathcal{M}$
{\it Ricci explodes} if there exists a geodesic $\gamma$ defined
for $t\in(0,T)$ so
$\lim_{t\rightarrow0}|\rho(\dot\gamma(t),\dot\gamma(t))|=\infty$.

\begin{example}\label{exm-6.3}\rm
For $1\le i\le 6$, let $\mathcal{N}_{i,\pm}:=\mathcal{M}_{f_{i,\pm}}$ where
$$\begin{array}{llll}
f_{1,-}(y)=-e^{-y},&f_{2,-}(y)=-e^{-y}+y,&f_{3,-}(y)=-e^{-y}-e^{-2y},\\
f_{1,+}(y)=e^y,&f_{2,+}(y)=e^y+y,&f_{3,+}(y)=e^y+e^{2y}\,.\end{array}$$
\begin{enumerate}
\item $\mathcal{N}_{1,-}$ is locally
homogeneous and Ricci explodes.
\item $\mathcal{N}_{2,-}$ is $1$-curvature homogeneous, not
$2$-curvature homogeneous, and  Ricci explodes.
\item $\mathcal{N}_{3,-}$ is $0$-curvature homogeneous, not $1$-curvature homogeneous, and Ricci explodes.
\item $\mathcal{N}_{1,+}$ is geodesically complete, and  homogeneous.
\item $\mathcal{N}_{2,+}$ is $1$-curvature homogeneous, not
$2$-curvature homogeneous,  and geodesically complete.
\item $\mathcal{N}_{3,+}$ is $0$-curvature homogeneous, not $1$-curvature homogeneous,  and geodesically complete.
\end{enumerate}\end{example}

\section*{Acknowledgments}
The research of M. Brozos-V\'azquez was partially supported by
project BFM 2003-02949 (Spain). The research of both M. Brozos-V\'azquez and  P. Gilkey was
partially supported by the Max Planck Institute for
Mathematics in the Sciences (Leipzig, Germany).

\end{document}